\documentclass[11pt]{article}
\usepackage{amssymb}
\usepackage{mathrsfs}
\def\flex#1{\mathrel{\mathop{\kern 0pt\hbox to 
10mm{\rightarrowfill}}\limits_{#1
\rightarrow \infty}}}

\def\Flex#1{\mathrel{\mathop{\kern 0pt\hbox to 
10mm{\rightarrowfill}}\limits_{#1
\rightarrow \infty}}}

\def\u{{\bf u}}
\def\v{{\bf v}}
\def\1{1\hspace{-1.2mm}\mbox{{\normalsize I}}}

\def\Dem{\vskip 2.5mm\noindent { Proof{}.}  }
\newtheorem{Th}{Theorem}

\newtheorem{Lemme}[Th]{Lemma}
\newtheorem{Cor}[Th]{Corollary}

\newtheorem{Prop}[Th]{Proposition}
\newcommand{\R}{\mathbb{R}}

\title{Oblique repulsion in the nonnegative quadrant}
\author{Dominique L\'epingle\\
   Universit\'e d'Orl\'eans, FDP-MAPMO}

\begin{document}
\maketitle

\begin{abstract}
We consider the differential system $\dot{x}=\alpha/x+\beta/y$, $\dot{y}=\gamma/x+\delta/y$ in the nonnegative quadrant. 
Here $\alpha$ and $\delta$ are positive, $\beta$ and $\gamma$ are real constants. Under some condition on the constants there exists
a unique global solution. The main difficulty is to prove uniqueness when starting at the corner of the quadrant.
\end{abstract} 

\section{Introduction.}

We are interested in the question of existence and uniqueness of the solution $\u(.)=(x(.),y(.))$ to the following integral system
\begin{equation}
 \label{eq:syst}
  \begin{array}{lll}
    x(t) & = & x+\alpha\int_0^t\frac{ds}{x(s)} + \beta \int_0^t\frac{ds}{y(s)} \\
    y(t) & = & y + \gamma\int_0^t \frac{ds}{x(s)} + \delta \int_0^t\frac{ds}{y(s)}
  \end{array}
\end{equation}
where $x(.)$ and $y(.)$ are continuous functions from $[0,\infty)$ to $[0,\infty)$ with the conditions 
\begin{equation}
  \begin{array}{lll}
   \int_0^t \1_{\{x(s)=0\}}ds = 0 & \qquad & \int_0^t \1_{\{y(s)=0\}}ds = 0 \\
   \int_0^t \1_{\{x(s)>0\}} \frac{ds}{x(s)} <\infty & \qquad  &  \int_0^t \1_{\{y(s)>0\}} \frac{ds}{y(s)} <\infty 
  \end{array}
\end{equation}
for any $t\geq 0$. Here $\alpha$, $\beta$, $\gamma$ and $\delta$ are four real constants with $\alpha>0$ and $\delta>0$.

 The system has a single singularity at each side of the nonnegative quadrant $S=\{(x,y): x\geq 0,y\geq 0 \}$ and a double singularity
at the corner ${\bf 0}=(0,0)$. We write $S^{\circ}:=S\setminus\{ {\bf 0}\}$ for the punctured quadrant.

We will note $\dot{x}(t)$ the derivative $dx(t)/dt$. So the integral system (\ref{eq:syst}) may be written as an initial-value problem
\begin{equation}
\label{eq:valini}
\begin{array}{lll}
  \dot{x} & = & \frac{\alpha}{x} + \frac{\beta}{y} \\
  \dot{y} & = & \frac{\gamma}{x} + \frac{\delta}{y} 
  \end{array}
\end{equation}
with the inital condition $(x(0),y(0))\in S$.

We first remark that if $\beta<0$, $\gamma<0$ and $\alpha \delta < \beta \gamma$, there exist $\lambda>0$ and $\mu>0$ such that
$\lambda \alpha + \mu \gamma <0$ and $\lambda \beta + \mu \delta <0$. Thus $z(t):=\lambda x(t) + \mu y(t)$ is decreasing, $\min{(x(t),y(t))} \rightarrow 0$  
and $\dot{z}(t) \rightarrow -\infty$ as $t\rightarrow t_f$ where $t_f<\infty$ and there is no solution.
 If $\beta <0$, $\gamma<0$ and $\alpha \delta = \beta \gamma$, then $v(t):= \alpha y(t) - \gamma x(t)$ remains equal to $\alpha y - \gamma x$ and there is a unique solution $(x(t),y(t))$ that converges to $(\frac{\gamma x - \alpha y}{\beta + \gamma}, \frac{\beta y - \delta x}{\beta + \gamma})$
 except if $(x,y)= {\bf 0}$ in which case there is no solution. 

From now on we will make the following hypothesis:
\[
  (H) \qquad \qquad \qquad \max{(\beta, \gamma) } \geq 0  \;\; \mbox{or} \;\; \beta \gamma< \alpha \delta .
\]
This is equivalent to the existence of $\lambda\geq 0$ and $\mu\geq 0$ such that $\lambda \alpha + \mu \gamma >0$
 and $\lambda \beta + \mu \delta > 0$. This last formulation amounts to saying that the matrix 
\[
   \begin {array}{lll}
     A & = & \left( \begin{array}{cc}
                              \alpha & \beta \\
                             \gamma & \delta
                            \end{array}
                   \right)
   \end{array}
\] 
is an $S$-matrix in the terminology of \cite{FP}. In the sequel, we fix a pair $(\lambda, \mu)$ with $\lambda>0$, $\mu>0$, such that 
$\lambda \alpha + \mu \gamma >0$ and $\lambda \beta + \mu \delta >0$.

The aim of this note is to prove the following result.

\begin{Th}
\label{Th:objet}
 Under condition $(H)$, there exists a unique solution to {\rm(\ref{eq:syst})} for any starting point $(x,y)\in S$.
\end{Th}

\section{Some preliminary lemmata.}

We begin with a comparison lemma.

\begin{Lemme}
\label{Th:comp}
 Let $x_1$ and $x_2$ be nonnegative continuous functions on $[0,\infty)$ which are solutions to the system
\[
  \begin{array}{lll}
    x_1(t) & = & v_1(t) + \alpha \int_0^t\frac{ds}{x_1(s)} \\
    x_2(t) & = & v_2(t) + \alpha \int_0^t\frac{ds}{x_2(s)}
  \end{array}
\]
where $\alpha>0$, $v_1$ and $v_2$ are continuous functions such that $0\leq v_1(0)\leq v_2(0)$, and $v_2-v_1$ is nondecreasing.
Then $x_1 \leq x_2$ on $[0,\infty)$.
\end{Lemme}  
\Dem
Assume there exists $t>0$ such that $x_1(t)>x_2(t)$. Set
\[
  \tau := \max{\{s\leq t: x_1(s)\leq x_2(s)\}}.
\]
Then
\[
  \begin{array}{ll}
  x_2(t)-x_1(t)& =  x_2(\tau)-x_1(\tau)+(v_2(t)-v_1(t))-(v_2(\tau)-v_1(\tau))+\alpha \int_{\tau}^t(\frac{1}{x_2(s)}-\frac{1}{x_1(s)})ds \\
                &\geq   0 ,
  \end{array}
\]
a contradiction. $\hfill  \blacksquare $

\begin{Lemme}
\label{Th:sysmixt}
Let the system
\begin{equation}
 \label{eq:sysmix} 
 \begin{array}{lll}
  \dot{x} & = & \frac{\alpha}{x}+ \phi(x,z) \\
  \dot{z} & = & \psi(x,z)
  \end{array}
\end{equation}
with $x(0)=x_0\geq 0$, $z(0)=z_0\in \R$, $\alpha>0$, $\phi$ and $\psi$ two Lipschitz functions on $\R_+\times \R$, 
and $| \phi| \leq c$ for some $c<\infty$. Then there exists a unique solution to {\rm(\ref{eq:sysmix})}. Moreover, for this solution, 
$x(t)>0$ for any $t>0$.
\end{Lemme}
\Dem
Assume first $x_0>0$. Then the system (\ref{eq:sysmix}) is Lipschitz on $[\min{\{x_0, \frac{\alpha}{c}\}}, \infty)\times \R$ 
and the solution does not step out of
this domain, so there is a unique global solution. When $x_0=0$, we let  $w_0(t)=0$, $ z_0(t)=z_0$ and for $n\geq 1$
\[
  \begin{array}{lll}
     w_n(t) & = & 2\alpha t + 2 \int_0^t \sqrt{w_{n-1}(s)}\, \phi(\sqrt{w_{n-1}(s)}, z_{n-1}(s)) ds \\
     z_n(t) & = & z_0 + \int_0^t \psi(\sqrt{w_{n-1}(s)},z_{n-1}(s)) ds .
  \end{array}
\]
Let $M>0$ and assume $|w_{n-1}(t)|\leq M$ on some interval $[0,T]$. Then, for $0\leq t\leq T$,
\[
  |w_n(t)| \leq T(2\alpha + 2c \sqrt{M})
\]
and this is again $\leq M$ for $T $ small enough. We also have $|z_n(t)|\leq M^{\prime}$ for any $n\geq 0$ for $T$ small enough.
Equicontinuity of $(w_n,z_n)_{n\geq 0}$ is easily verified and from the Arzela-Ascoli theorem it follows there exists a subsequence 
$(w_{n_k},z_{n_k})$ converging on $[0,T]$ to a solution $(w,z)$ of the system
\begin{equation}
  \label{eq::syscar}
   \begin{array}{lll}
    \dot{w} & = & 2\alpha + 2 \sqrt{w}\, \phi(\sqrt{w},z) \\
    \dot{z} & = & \psi(\sqrt{w},z) 
   \end{array}
\end{equation} 
with the initial conditions $w(0)=0$, $z(0)=z_0$. For small $T$, $\dot{w}>0$ on $[0,T]$. Set now $x(t)=\sqrt{w(t)}$. Then $(x,z)$  is a solution to
(\ref{eq:sysmix}) on $[0,T]$ with $x(T)>0$. We may extend the solution to $[0,\infty)$ by using the above result with $x_0>0$. 

We now prove uniqueness. Let $(x,z)$ and $(x^{\prime},z^{\prime})$ be two solutions of (\ref{eq:sysmix}). Then
\[
  \begin{array}{ll}
     & (x(t)-x^{\prime}(t))^2 + (z(t)-z^{\prime}(t))^2 \\
  = & 2 \alpha \int_0^t(x(s)-x^{\prime}(s))(\frac{1}{x(s)}-\frac{1}{x^{\prime}(s)})ds +2 \int_0^t(x(s)-x^{\prime}(s))
   (\phi(x(s),z(s))-\phi(x^{\prime}(s),z^{\prime}(s))) ds \\
   & + 2 \int_0^t (z(s)-z^{\prime}(s))(\psi(x(s),z(s))-\psi(x^{\prime}(s),z^{\prime}(s)))ds \\
  \leq & 4 L  \int_0^t((x(s)-x^{\prime}(s))^2 + (z(s)-z^{\prime}(s))^2) ds
  \end{array}
\]
where $L$ is the Lipschitz constant of $\phi$ and $\psi$. Uniqueness follows from Gronwall's inequality. $\hfill  \blacksquare$

\begin{Lemme}
\label{Th:crois}
Let $\u(.)=(x(.),y(.))$ be a solution to (\ref{eq:syst}) and let $\nu=(\lambda,\mu)$. Then the  function
$z(t):=\nu. \u(t)=\lambda x(t) + \mu y(t)$ is increasing on $[0,\infty)$ and 
we have $\u(t)\in S^{\circ}$ for any $t>0$.
\end{Lemme}
\Dem
Recall that condition $(H)$ is in force. We easily check that $\dot{z}(t)$ is positive. $\hfill \blacksquare$     

\section{Existence. Case $x=0,y=0.$}

There is an explicit solution to (\ref{eq:syst}) when the starting point is the corner.

\begin{Prop}
There is a solution to (\ref{eq:syst}) with initial condition ${\bf 0}$ given by 
\begin{equation}
 \label{eq:coin}
  \begin{array}{lll}
   x(t) & = & c\sqrt{t} \\
   y(t) & = & d\sqrt{t}
  \end{array}
\end{equation}
where
\begin{equation}
 \label{eq:const}
 \begin{array}{lll}
    c & = & (2\alpha + \frac{\beta}{\delta} (\beta-\gamma + \sqrt{(\beta-\gamma)^2 + 4 \alpha \delta}))^{1/2} \\
    d & = & (2\delta + \frac{\gamma}{\alpha} (\gamma-\beta + \sqrt{(\beta-\gamma)^2 + 4 \alpha \delta}))^{1/2}.
 \end{array}
\end{equation}
\end{Prop}
\Dem
Writing down $x(t)=c\sqrt{t}$ and $y(t)=d\sqrt{t}$ we have to solve the system
\[
  \begin{array}{lll}
    \frac{c}{2} & = & \frac{\alpha}{c} + \frac{\beta}{d} \\
    \frac{d}{2} & = & \frac{\gamma}{c} + \frac{\delta}{d} ;
  \end{array}
\]
We first compute
\begin{equation}
 \label{eq:valeur}
 \frac{d}{c} = \frac{\gamma-\beta + \sqrt{(\beta-\gamma)^2 + 4 \alpha \delta}}{2 \alpha}
\end{equation}
and then obtain (\ref{eq:const}) provided that 
\[
  \begin{array}{lll}
  C & = & 2\alpha + \frac{\beta}{\delta} (\beta-\gamma + \sqrt{(\beta-\gamma)^2 + 4 \alpha \delta}) \\
  D & = & 2\delta + \frac{\gamma}{\alpha} (\gamma-\beta + \sqrt{(\beta-\gamma)^2 + 4 \alpha \delta})
  \end{array}
\]
are positive. If $\beta\geq 0$, $C$ is clearly positive. This is also true if $\beta<0$ and $\beta \gamma < \alpha \delta$ since $C$
may be written
\[
   C = \frac{4 \alpha (\alpha \delta -\beta \gamma)}
                          {2 \alpha \delta- \beta \gamma + \beta^2 - \beta\sqrt{4 (\alpha \delta - \beta \gamma) + (\beta + \gamma)^2}} . 
\]
The proof for $D$ is similar. $\hfill \blacksquare$

Uniqueness in this case is more involved and will be treated in the last section. We only remark for the moment that the system (\ref{eq:valini})
with $\alpha=\delta=0$, $\beta>0$, $\gamma>0$ and initial value ${\bf 0}$ has a one-parameter family of solutions.

\section{Angular behavior.}

We are now in a position to study the behavior of $\frac{y(t)}{x(t)}$. For any $\u=(x,y)\in S^{\circ}$ we set
\[
  \theta(\u)= \arctan {\frac{y}{x}}.
\]
We also set
\[
 \u_{\ast}=(x_{\ast},y_{\ast}):= \left(\frac{c}{\lambda  c+\mu  d},\frac{d}{\lambda  c+ \mu d}\right).
\]

\begin{Prop}
\label{Th:arctan}
Let $\u(.)$ be a solution to  (\ref{eq:syst}) starting at $\u=(x,y)\in S^{\circ}$. Then for any $t>0$
\begin{enumerate}
 \item
\[
 \begin{array}{lllll}
  \frac{d\theta(\u(t))}{dt} >0 &\;\;{\rm and} &\;\; \theta(\u(t))<\theta(\u_{\ast}) & \qquad {\rm if} & \quad \theta(\u)<\theta(\u_{\ast})  \\
   \frac{d\theta(\u(t))}{dt} =0 &\;\;{\rm and} &\;\; \theta(\u(t))=\theta(\u_{\ast}) & \qquad {\rm if} & \quad \theta(\u)=\theta(\u_{\ast})  \\
   \frac{d\theta(\u(t))}{dt} <0 &\;\;{\rm and} &\;\; \theta(\u(t))>\theta(\u_{\ast}) & \qquad {\rm if} & \quad \theta(\u)>\theta(\u_{\ast})  .
 \end{array}
\]
\item
 \begin{equation}
\label{eq:mino}
\begin{array}{lll}
  x(t) & \geq & \min{\left(x,c \frac{\lambda x + \mu y}{\lambda c+ \mu d}\right)}  \\
  y(t) & \geq & \min{\left(y,d \frac{\lambda x + \mu y}{\lambda c+ \mu d}\right)} .
\end{array}
\end{equation}
\end{enumerate}
\end{Prop}
\Dem
From Lemma \ref{Th:crois} we know that $\u(t)\in S^{\circ}$ for any $t\geq 0$.
\begin{enumerate}
 \item We compute
\begin{equation}
  \label{eq:angul}
   \frac{d\theta(\u(t))}{dt} = \frac{1}{x^2(t)+y^2(t)}\left(\frac{d}{c}-\frac{y(t)}{x(t)}\right)\left[\alpha+\frac{x(t)(\beta-\gamma+
 \sqrt{(\beta-\gamma)^2 + 4 \alpha \delta})}{2 y(t)}\right] 
\end{equation}
and the conclusion follows.
 \item
 Let ${\bf a},{\bf b}\in S^{\circ}$ with $0\leq\theta({\bf a})<\theta(\u_{\ast})$,
 $\theta(\u_{\ast})<\theta({\bf b})\leq\frac{\pi}{2}$ and let $l>0$. We set 
\begin{equation}
\label{eq:parti}
  \begin{array}{lll}
  A & = & \{ \v \in S^{\circ}: \theta({\bf a }) \leq \theta(\v) \leq \theta(\u_{\ast}) \} \\
  B & = & \{ \v \in S^{\circ}: \theta({\bf b }) \geq \theta(\v) \geq \theta(\u_{\ast}) \} .
  \end{array}
\end{equation}
It follows from above that any solution starting from $A$ stays in $A$, and the same is true for 
$B$.
  If $\u\in A$,
\[
  x(t)\geq -\frac{\mu}{\lambda}y(t) + (x+\frac{\mu}{\lambda}y) \geq -\frac{\mu d}{\lambda c}x(t)  + x +\frac{\mu}{\lambda}y
\]
and therefore 
\[ 
 x(t)\geq c \frac{\lambda x + \mu y}{ \lambda c + \mu d}.
\]
If  $\u\in B$,
\[
  x(t)\geq \frac{x}{y}y(t) \geq \frac{x}{y}(-\frac{\lambda}{\mu}x(t) + \frac{\lambda x+\mu y}{\mu})
\]
and therefore 
\[
  x(t)\geq x .
 \]
Same estimations for $y(t)$. $\hfill \blacksquare$ 
\end{enumerate}

\begin{Cor}
 \label{Th:pente}
Let $\u(.)$ be a solution to  (\ref{eq:syst}). Then
\[
  \lim_{t\rightarrow \infty} \theta(\u(t)) = \theta(\u_{\ast}),
 \qquad
i.e.
\qquad
  \lim_{t\rightarrow \infty} \frac{y(t)}{x(t)} = \frac{d}{c} .
\]
\end{Cor}
\Dem
If $\u=(x,y) \in S^{\circ}$, this is an easy consequence of (\ref{eq:angul}). If $\u=(0,0)$ we may apply Lemma \ref{Th:crois} and then
(\ref{eq:angul}). $\hfill \blacksquare $

\section{Existence and uniqueness. Case $x>0,y>0$. }

\begin{Prop}
\label{Th:int}
There exists a unique solution $\u(.)$  to  (\ref{eq:syst}) starting at $\u=(x,y)$ with $x>0,y>0$. It satisfies 
$x(t)>0,y(t)>0$ for any $t\geq 0$.
\end{Prop}
\Dem
We now assume $\theta({\bf a })>0$ and $\theta({\bf b})< \frac{\pi}{2}$ in (\ref{eq:parti}).
 Let $l>0$ and ${\bf \nu}=(\lambda,\mu)$.  We set
$ L  :=  \{ \v \in S^{\circ}: {\bf \nu}.\v \geq l\}$.
From Lemma \ref{Th:crois} and Proposition \ref{Th:arctan} we know that any solution starting from $A\cap L$ stays in $A \cap L$, and the same is true for 
$B\cap L$. As the system is Lipschitz in $A\cap L$ and in $B\cap L$, there is a unique global solution to (\ref{eq:syst}) in both cases. $\hfill \blacksquare$ 
 
\section{Existence and uniqueness. Case $x=0,y>0$.}

\begin{Prop}
\label{Th:mixte} 
There exists a unique solution $\u(.)$  to  (\ref{eq:syst}) starting at $\u=(x,y)$ with $x=0,y>0$. It satisfies $x(t)>0, y(t)>0$ for
any $t>0$.
\end{Prop}
\Dem
Let $\varepsilon\in(0,y\frac{\mu d}{\lambda c+ \mu d})$. We define on $\R_+\times \R$
\[
  \psi_{\varepsilon}(x,z) := \frac{1}{\max{(\gamma x +z,\alpha  \varepsilon)}} .
\]
We apply Lemma \ref{Th:sysmixt} to obtain a unique solution $x_{\varepsilon}(.), z_{\varepsilon}(.)$ to
\begin{equation}
  \begin{array}{lll}
    x_{\varepsilon}(t) & = & \alpha \int_0^t\frac{ds}{x_{\varepsilon}(s)} + \alpha^m \beta \int_0^t\psi_{\varepsilon}(x_{\varepsilon}(s), z_{\varepsilon}(s))ds \\
    z_{\varepsilon}(t) & = & \alpha y + \alpha(\alpha \delta-\beta \gamma)  \int_0^t\psi_{\varepsilon}(x_{\varepsilon}(s), z_{\varepsilon}(s))ds .
  \end{array}
\end{equation}
Let
\[
  \begin{array}{lll}
  y_{\varepsilon}(t) & = & \frac{1}{\alpha}(\gamma x_{\varepsilon}(t) + z_{\varepsilon}(t)) \\
  \tau_y(\varepsilon) & = & \inf{\{t>0:y_{\varepsilon}(t)<\varepsilon\}} .
  \end{array}
\]
On the interval $[0,\tau_y(\varepsilon)]$,  $(x_{\varepsilon}(.), y_{\varepsilon}(.))$ is the unique
 solution to (\ref{eq:syst}). From (\ref{eq:mino}) 
 we know that $y_{\varepsilon}(t)>\varepsilon$ on this interval. Thus  $\tau_y(\varepsilon))=\infty$
and $(x(.),y(.)):=(x_{\varepsilon}(.), y_{\varepsilon}(.))$ is the unique global solution to  (\ref{eq:syst}). $\hfill \blacksquare$

\section{Path behavior.}

  Let us note $\u(t,\u_0)$ the solution to (\ref{eq:syst}) starting 
at $\u_0\in S^{\circ}$. Using Gronwall's inequality as in the proof of uniqueness, it is easily seen that for any $t>0$ 
the solution 
$\u(t,\u_0)$ 
continuously depends
on the initial condition $\u_0$. It has the 
Scaling Property:
\[
   {\rm (SC)}  \qquad  \qquad \u(r^2t,\u_0) = r\u(t,\frac{\u_0}{r})
\]
for any $r>0$, $t\geq 0$, $\u_0\in S^{\circ}$. Using Lemma \ref{Th:crois} we also note that any  solution $\u(.)$ to   (\ref{eq:syst}) has the Semi-group Property:
\[
   {\rm (SG)}  \qquad  \qquad \u(s+t)=\u(t,\u(s))
\]
for any $s>0$ and $t\geq 0$. With Proposition \ref{Th:int} and Proposition \ref{Th:mixte} this entails that $x(t)>0$ and $y(t)>0$ 
for any $t>0$.  We now set for any 
$r>0$:
\[
  \begin{array}{lll}
  L_r & := & \{ \u=(x,y): x>0,y>0, \nu.\u=r\} \\
 \overline{L_r} & := & \{ \u=(x,y): x\geq 0,y\geq 0, \nu.\u=r\}
  \end{array}
\]

\begin{Lemme}
 \label{Th:attein}
  Let $\u(.)$  be a solution to (\ref{eq:syst}) starting at $\u_0\in S$ with $\nu.\u_0 \leq r$. We set
\[
  \tau(r) := \inf{\{t\geq 0 : \nu.\u(t)= r\}} .
\]
Then 
\[
  \tau(r) \leq \frac{r^2}{2[\lambda(\lambda\alpha+\mu \gamma)+ \mu (\lambda \beta + \mu \delta)]} .
\]
\end{Lemme}

\Dem
Set
\[
  z(t):= \nu.\u(t) .
\]
As
\[
  z(t) = \nu.\u_0 + [\lambda(\lambda\alpha+\mu \gamma)+ \mu (\lambda \beta + \mu \delta)] \int_0^t \frac{ds}{z(s)} + \int_0^t f(s) ds
\]
with
\[
  f(s)=  \mu (\lambda\alpha+\mu \gamma) \frac{y(s)}{x(s)}+ \lambda (\lambda \beta + \mu \delta) \frac{x(s)}{y(s)} >0
\]
for $s>0$, it follows from Lemma \ref{Th:comp} that $z(t)\geq w(t)$ 
where
\[
  w(t)= \nu.\u_0 +  [\lambda(\lambda\alpha+\mu \gamma)+ \mu (\lambda \beta + \mu \delta)] \int_0^t\frac{ds}{w(s)},
\]
and then 
\[
  z^2(t) \geq w^2(t) = 2  [\lambda(\lambda\alpha+\mu \gamma)+ \mu (\lambda \beta + \mu \delta)] t + (\nu.\u_0)^2 .
\]
The conclusion follows. $\hfill \blacksquare$

We now define $q: \overline{L_1}\rightarrow L_1$ by
\[
  q(\u_1) = \frac{1}{2} \u(\tau(2),\u_1)
\]
where
\begin{equation}
 \label{eq:to}
  \tau(2) = \inf{\{t\geq 0: \nu.\u(t,\u_1)=2\}}
\end{equation}
is finite from the above Lemma. Let now $r>0$ and $\u \in\overline{L_r}$. From (SC), the geometric paths in $S$ of $\u(.,\u)$ and $r \u(.,\frac{\u}{r})$ are identical. 
Therefore
\[
   \u(\tau(2r),\u) = r \u(\tau(2),\frac{\u}{r})
\]
where in this equality $\tau(2r)$ is relative to $\u(.,\u)$ and $\tau(2)$ is relative to $\u(.,\frac{\u}{r})$. Thus
\[
  q(\frac{\u}{r}) = \frac{1}{2 r} \u(\tau(2r),\u).
\]
Iterating and using (SG), we get for any $n\geq 1$
\begin{equation}
 \label{eq:iter}
 q^n(\frac{\u}{r}) = \frac{1}{2^n r}\u(\tau(2^n r),\u) .
\end{equation}

\begin{Prop}
\label{Th:conv}
There exists $k\in (0,1)$ such that for any $\u_1\in \overline{L_1}$ 
\[
  |q(\u_1)-\u_{\ast}| \leq k\, |\u_1-\u_{\ast}|.
\]
\end{Prop}
\Dem
From Proposition \ref{Th:arctan} we know that $q$ has a unique invariant point $u_{\ast}$.
We consider the solution $\u(t,\u_1)=(x(t),y(t))$ on the time interval $[0,\tau(2)]$, where $\tau(2)$ was defined in (\ref{eq:to}). We first assume that 
\[
  \frac{y_1}{x_1}< \frac{y_{\ast}}{x_{\ast}}= \frac{d}{c} .
\]
We note for further use that
\[
  \begin{array}{lllllll}
  x_{\ast} & < & x_1 & \leq & \frac{1}{\lambda} & & \\
  0 & \leq & y_1 & < & y_{\ast} & < & \frac{1}{\mu} .
  \end{array} 
\] 
We set
\[
  \u_2 = (x_2,y_2) :=\u_1 +\frac{(\alpha y_1 + \beta x_1, \gamma y_1 + \delta x_1)}{\lambda(\alpha y_1 + \beta x_1) + \mu (\gamma y_1 + \delta x_1)} .
\]{
Then, $2 \u_{\ast}$, $\u_{\ast}+\u_1$ and $\u_2 \,\in L_2$. Setting for $z\in [0,\infty]$
\[
  h(z) = \frac{\gamma z + \delta}{\alpha z + \beta}
\]
we compute
\begin{equation}
 \label{eq:deriv}
  \frac{dh}{dz}(z) = \frac{\beta \gamma - \alpha \delta}{(\alpha z + \beta)^2}\,.
\end{equation}
From Proposition \ref{Th:arctan} we know that for any $t\in [0,\tau(2)]$
\begin{equation}
\label{eq:inega}
  \frac{y(t)}{x(t)} < \frac{d}{c}.
\end{equation}
When $\alpha \delta > \beta \gamma$, it follows from (\ref{eq:deriv}) and (\ref{eq:inega}) that
\[
  \frac{\dot{y}(t)}{\dot{x}(t)} = h(\frac{y(t)}{x(t)}) > h(\frac{d}{c}) = \frac{d}{c}
\]
and then $2 q(\u_1)$ belongs to the open interval $(2 \u_{\ast}, \u_{\ast} + \u_1)$ on $L_2$. 
Therefore,
\[
  |2 q(\u_1)- 2 \u_{\ast}| < |\u_1- \u_{\ast}|.
\]
When $\alpha \delta = \beta \gamma$, the path of the solution is a straight half-line with slope $\frac{d}{c}$ and 
\[
  |2 q(\u_1) -2 \u_{\ast}|  = |\u_1-\u_{\ast}| .
\] 
When $\alpha \delta < \beta \gamma$,  $\frac{\dot{y}(t)}{\dot{x}(t)}$  is increasing on $[0, \tau(2)]$ and then
\[
  \frac{\gamma y_1 + \delta x_1}{\alpha y_1 + \beta x_1} \leq \frac{\dot{y}(t)}{\dot{x}(t)} < \frac{d}{c} .
\]
As a result, $2 q(\u_1)$ belongs to the open interval $(\u_{\ast} + \u_1, \u_2)$ on $L_2$. Moreover, using the relation $\lambda x_1 + \mu y_1 = 1$ twice, we get
\[
  \begin{array}{lll}
   2 x_1-x(\tau(2)) & > & 2 x_1 - x_2 \\
        & = & x_1 - \frac{\alpha y_1 + \beta x_1}{\lambda(\alpha y_1 + \beta x_1) + \mu (\gamma y_1 + \delta x_1)} \\
        & = & \frac{\alpha \lambda x_1 y_1 + \beta \lambda x_1^2 + \gamma \mu x_1 y_1 + \delta \mu x_1^2 -\alpha y_1 -\beta x_1 }
                  {\lambda(\alpha y_1 + \beta x_1) + \mu (\gamma y_1 + \delta x_1)} \\
        & = & \mu \frac{-\alpha y_1^2 + (\gamma-\beta )x_1 y_1 + \delta x_1^2}{\lambda(\alpha y_1 + \beta x_1) + \mu (\gamma y_1 + \delta x_1)} \\
       & = & \frac{\alpha \mu}{\lambda(\alpha y_1 + \beta x_1) + \mu (\gamma y_1 + \delta x_1)} (x_1 \frac{d}{c}-y_1)
       \left (y_1+ \frac{\beta - \gamma + \sqrt{(\beta-\gamma)^2+ 4 \alpha \delta}}{2 \alpha} x_1\right) \\
       & = & \frac{\alpha (\lambda c + \mu d)}{c[\lambda(\alpha y_1 + \beta x_1) + \mu (\gamma y_1 + \delta x_1)]}
       (x_1 - x_{\ast})\left(y_1+ \frac{\beta - \gamma + \sqrt{(\beta-\gamma)^2+ 4 \alpha \delta}}{2 \alpha} x_1\right) . 
\end{array}
\]
In the same way,
\[
  \begin{array}{lll}
 y(\tau(2))-2 y_1 & > & y_2-2 y_1 \\
        & = & \frac{\alpha (\lambda c + \mu d)}{c[\lambda(\alpha y_1 + \beta x_1) + \mu (\gamma y_1 + \delta x_1)]}
       ( y_{\ast}-y_1) \left(y_1+ \frac{\beta - \gamma + \sqrt{(\beta-\gamma)^2+ 4 \alpha \delta}}{2 \alpha} x_1\right) . 
  \end{array}
\]
Setting
\[
  \begin{array}{lllll}
  k_1 & = &  \frac{\lambda \mu (\beta - \gamma + \sqrt{(\beta-\gamma)^2+ 4 \alpha \delta})}
{4[\lambda(\lambda \alpha  + \mu \gamma) + \mu (\lambda \beta + \mu\delta )]} & > & 0
  \end{array}
\]
we obtain
\[
  |2 q(\u_1)-2 \u_1 | > 2k_1\,|\u_1 - \u_{\ast}|
\]
and then
\[
  \begin{array}{lll} 
     |q(\u_1)-\u_{\ast}| & = & |\u_1-\u_{\ast}| -|q(\u_1)-\u_1| \\
   & < & (1-k_1)\, |\u_1-\u_{\ast}| .
   \end{array}
 \]
If now $\frac{y_1}{x_1}> \frac{d}{c}$, in the same way there exists $k_2>0$ such that
\[
   |q(\u_1)-\u_{\ast}| < (1 - k_2)\,|\u_1-\u_{\ast}|
\]
We may take $k= 1-\min(k_1,k_2)          \qquad \qquad \hfill \blacksquare $

\section{Uniqueness. Case $x=0,y=0$.}

Existence  was proven in Section 3. We may now conclude the proof of Theorem \ref{Th:objet}.

\begin{Prop}
The solution given by (\ref{eq:coin}) is the unique solution to  (\ref{eq:syst}) starting at ${\bf 0}$. 
\end{Prop}

\Dem
 Let $\u(.)$ be a solution to (\ref{eq:syst}) starting at ${\bf 0}$. 
For any $n\geq 1$ and $s>0$,
\[
  \u(\tau(s))=\u(\tau(s),\u(\tau(s2^{-n})))
\]
where $\tau(s)$ in the l.h.s. is relative to $\u(.)$ and  $\tau(s)$ in the r.h.s. is  relative to $\u(.,\u(\tau(s2^{-n})))$. 
We may apply (\ref{eq:iter}) with
$r=s2^{-n}$ and $\u=\u(\tau(s2^{-n}))$. We obtain
\[
  \frac{\u(\tau(s))}{s} = q^n(\frac{\u(\tau(s2^{-n}))}{s2^{-n}}) .
\]
From Proposition \ref{Th:conv} (or directly from (\ref{eq:angul}) it follows that the r.h.s. converges to $\u_{\ast}$ as $n\rightarrow \infty$. Thus for any $s>0$
\[
   \frac{\u(\tau(s))}{s}=\u_{\ast}
\]
and this implies 
\[
  \frac{y(\tau(s))}{x(\tau(s))}= \frac{d}{c} .
\]
From Lemma \ref{Th:attein} we know that $\tau$ is one-to-one from $[0,\infty)$ to $[0,\infty)$ , and thus
for any $t>0$,
\[
  \frac{y(t)}{x(t)}= \frac{d}{c}.
\]
Going back to the system (\ref{eq:syst}) we conclude that 
\[
  \begin{array}{lll} 
   x(t) & = & c \sqrt{t}  \\
  y(t) & = & d \sqrt{t} . \qquad \qquad \hfill \blacksquare
  \end{array} 
\]

\end{document}